\documentclass[notitlepage,leqno,10pt]{amsart}
\usepackage{amsfonts,amssymb,cite}
\usepackage[latin1]{inputenc}
\usepackage{enumerate}
\usepackage{amssymb}
\catcode`\@=11 \@addtoreset{equation}{section}

\catcode`\@=12
\usepackage{latexsym}
\usepackage{textcomp}
\usepackage{amsmath}
\usepackage{amsfonts}
\usepackage{amssymb}
\usepackage{mathrsfs}

\newcommand{\D }{\Delta }

\newcommand{\e }{\varepsilon }

\renewcommand{\l }{\lambda }

\newcommand{\n }{\nabla }
\newcommand{\vp }{\varphi }

\newcommand{\s }{\sigma }

\renewcommand{\th }{\theta }
\renewcommand{\o }{\omega }
\renewcommand{\O }{\Omega }

\newcommand{\ov}{\overline}

\newcommand{\be}{\begin{equation}}
\newcommand{\ee}{\end{equation}}

\newcommand{\R}{\mathbb{R}}
\newcommand{\N}{\mathbb{N}}

\newcommand{\de}{\partial}

\newcommand{\ti}{\tilde}

\newcommand{\ra}{{\rangle}}
\newcommand{\la}{{\langle}}

\newcommand{\tz}{\tilde{z}}

\newcommand{\calD }{\mathcal{D}}

\newtheorem{Theorem}{Theorem}[section]
\newtheorem{Lemma}[Theorem]{Lemma}
\newtheorem{Proposition}[Theorem]{Proposition}

\def\proof{\noindent{{\bf Proof. }}}
\def\square{\vbox{
    \hrule height .4pt
    \hbox{\vrule width .4pt height 7pt \kern 7pt
       \vrule width .4pt}
    \hrule height .4pt }}

\def\QED{\hfill {$\square$}\goodbreak \medskip}

\def\R{{\mathbb R}}

\font\sc=cmcsc9  \textwidth=14truecm
\title{ The  role of the mean curvature in a  Hardy-Sobolev trace inequality}

\begin{document}

\author{Mouhamed Moustapha Fall, Ignace Aristide Minlend, El hadji Abdoulaye Thiam}
\address{\noindent African Institute for Mathematical Sciences (A.I.M.S.) of Senegal KM 2, Route de Joal, B.P. 1418 Mbour, S\'en\'egal }
\email{\small{mouhamed.m.fall@aims-senegal.org, ignace.a.minlend@aims-senegal.org, elhadji@aims-senegal.org} }

\date{}
\maketitle

\textbf{Abstract:} 
The Hardy-Sobolev trace inequality can be obtained via Harmonic extensions on the half-space of the Stein and Weiss weighted Hardy-Littlewood-Sobolev inequality.
In this paper we  consider a bounded domain and study the influence of the boundary mean curvature  in the  Hardy-Sobolev trace inequality
on the underlying domain. We  prove  existence of minimizers when the mean curvature is negative at  the singular point of the Hardy potential.\\


\section{Introduction}
The weighted Stein and Weiss  inequality (see \cite{StWe}) states, in particular, that there exists a constant $C(N,s)>0$ such that
$$
C(N,s) \displaystyle \left(\int_{  \R^{N}}|x|^{-s}|u|^{q(s)} \right)^{\frac{2}{q(s)}}dx
\leq   \int_{\R^N}|\xi| \widehat{u}^2d\xi \quad\forall u\in C^\infty_c(\R^N),
$$
where $ s\leq 1$,  $q(s)=\frac{2(N-s)}{N-1}$ and 
$$
\widehat{u}(\zeta)= \frac{1}{(2\pi)^{\frac{N}{2}}}\,
\int_{\R^N}e^{-\imath\xi\cdot{x} }u(x)dx
$$
is the Fourier transform of $u$.
We consider the Hardy-Sobolev trace constant which is given by 
\be\label{eq:Sofs}
S(s):=\inf_{u\in C^\infty_c(\R^N)}
\frac{ \displaystyle \int_{\R^N}|\xi| \widehat{u}^2d\xi}{ \displaystyle \left(\int_{  \R^{N}}|x|^{-s}|u|^{q(s)} \right)^{\frac{2}{q(s)}}
dx    }.
\ee
Denote by $$\R^{N+1}_+=\left\{z=(z^1,\ti{z})\in \R^{N+1}\quad:\quad z^1>0 \right\}$$  with boundary
 $\R^N\times\{0\}\equiv \R^N$.
We denote  and henceforth define $\calD:=\calD^{1,2}(\ov{\R^{N+1}_+})$ the completion of $C^\infty_c(\ov{\R^{N+1}_+} )$  with respect to the norm
$$
u \mapsto \int_{\R^{N+1}_+}|\n u|^2dz.
$$
Classical argument of harmonic extension, (see for instance \cite{CSilv} for generalizations) yields
\be\label{eq:CKNtrace} 
S(s)=
\inf_{u\in \calD}\frac{\displaystyle \int_{  \R^{N+1}_+}|\n u|^2 dz }
{\displaystyle \left(\int_{ \de \R^{N+1}_+}|\ti{z}|^{-s}|u|^{q(s)}\, d\ti{z} \right)^{\frac{2}{q(s)}}}.
\ee
Note that for $s=0$ then $q(0)=:2^{\sharp}$, the critical Sobolev exponent while  $S(0)$ coincides with the  Sobolev trace constant studied by Escobar \cite{Esc} and Beckner \cite{Beck} wiht applications in the Yamabe problem with prescribed mean curvature.
Existence of symmetric decreasing minimizers for the quotient  $S(s)$ in \eqref{eq:Sofs} were obtained by Lieb [\cite{Lieb}, Theorem 5.1]. We also quote the works \cite{NR1,NR2} for the existence of minimizers in critical Sobolev trace inequalities. If $s=1$, we recover $S(1)=2 \frac {\Gamma^2({\frac{N+1}{4}}) }{\Gamma^2(\frac{N-1}{4})}$, the relativistic Hardy constant (see e.g. \cite{Herbst}) which is never achieved in $\calD$. In this case, it is expected that there is no influence of the curvature in
comparison with the   works on Hardy inequalities with singularity at the boundary or in Riemannian manifolds, see \cite{Fall,Thiam}. \\
Let $\O$ be a smooth domain of $\R^{N+1}$, $N\geq 2$ with $0\in\de\O$. We consider  $(\de\O,\ti{g})$ as a  Riemaninan manifold, with
  Riemannian metric $\ti{g}$ induced by $\R^{N+1} $ on $\de \O$. Let  $d$ denote the  Riemannian distance in $(\de\O,\ti{g})$.
A classical  argument of partitioning of unity (see Lemma \ref{lem:HSt} below)  yields
the existence of  a constant $C(\O)>0$ such that
the following inequality
$$
C(\O){\left(\int_{\de\O}d^{-s}(\s) |u|^{q(s)}d\s\right)^{\frac{2}{q(s)}}} \leq \int_{\O}|\n u|^2dx
+\int_{\O}| u|^2dx \qquad \forall u\in H^1(\O).
$$

 Our aims in this paper is to study the existence of minimizers for the following quotient:

\be\label{eq:gOdef-i}
S(s,\O):= \inf_{u\in H^1(\O)}\frac{  \displaystyle  \int_{\O}|\n u|^2dx
+\int_{\O}| u|^2dx         }{  \displaystyle  \left(\int_{\de\O}d^{-s}(\s) |u|^{q(s)}d\s\right)^{\frac{2}{q(s)}}},
\ee
for $s\in[0,1)$.
%
Our main  result is the following:
\begin{Theorem}\label{th:exts-Om-i}
Let $\O$ be a smooth domain of $\R^{N+1}$, $N\geq 3$ with $0\in\de\O$ and $s\in[0,1)$. Assume that the mean curvature of $\de\O$ at 0 is negative.
Then   $S(s,\O)< S(s)$ and $S(s,\O) $ is achieved by a positive function $u\in H^1(\O)$ satisfying
\begin{align*}
\begin{cases}
\displaystyle-\D u+u=0&\qquad \textrm{ in } \O \vspace{3mm}\\
\displaystyle \frac{\de u}{\de \nu}=S(s,\O)\, d^{-s}(\s)\,  u^{q(s)-1}& \qquad \textrm{ on } \de \O,
\end{cases}
\end{align*}
where $\nu$ is the unit outer normal of $\de\O$.
\end{Theorem}

In the literature,  several authors studied the influence of curvature in the Hardy-Sobolev inequalities in Euclidean space and in Riemmanian manifolds, see  \cite{LiLin,GK,CL,GR1,GR2,GR3,GY,DN,Jaber} and the references there in. For instance, consider the Hardy-Sobolev constant:
\be\label{eq:CKN}
Q(s,\O):=\inf_{w\in H^1_0(\O)}\frac{\displaystyle \int_{\O}|\n w|^2dx   }{\displaystyle  \left(\int_{\O}|x|^{-s} |w|^{2(s)} dx \right)^{2/2(s)}},
\ee
with   $s\in(0,2)$ and $2(s)=\frac{2(N-s)}{N-2} $.
The role of the  local geometry  $\de \O$ at $0$ in the study of minimizers for $ Q(s,\O)$ was first  investigated  by Ghoussoub and Kang in \cite{GK}.
In \cite{GK} the authors showed that if all the principal curvatures of $\de \O$ at $0$ are negative then $Q(s,\O)< Q(s,\R^{N+1}_+) $ and it is achieved. This result were improved by Ghoussoub and Roberts assuming only that the mean curvature is negative at 0 while $N\geq 4$. 
Later Demyanov and Nazarov in \cite{DN} constructed domains in which $ Q(s,\O)$ is achieved wile the mean curvature of $\de \O$ 
at $0$ is not negative. Actually,   by the results in \cite{DN}, extremals for $Q(s, \O)$ exists if $\O$ is "average concave in a neighborhood of the origin". Later on, in the same year,    Ghoussoub and Robert \cite{GR1,GR2} used   refined blow-up analysis to prove existence of an
extremal for $Q(s, \O)$  provided the mean curvature of $\de\O$ is negative at 0.

 Recently Chern and Lin in \cite{CL}
proved that if the mean curvature of $\de \O$ at $0$ is negative then $Q(s,\O)< Q(s,\R^{N+1}_+) $ for $N\geq 2$ 
 and in these cases, $ Q(s,\O)$ is attained. See also the recent work of  Li and Lin \cite{LiLin} for generalizations.\\

We point out that the study of the effect of the curvature in the Hardy-Sobolev trace inequality seems to be quite rare in the literature while the Sobolev trace  ($s=0$) inequality have been  intensively studied in the last years, see for instance \cite{NR1,NR2}. According  to the authors level of information, the paper is one of the first dealing with this question. We would like to emphasize that our argument of proof (based on blow up analysis, in Proposition \ref{prop:exist_smO}) is different from those in the papers cited above.  The main observation is that, dealing with "pure" Hardy-Sobolev mnimization problem $s\in(0,1]$, one can depict a sequence of radii $r_n\to 0$ where, if  blow up occur, then concentration can only happen in $B(0,r_n)\cap \de \O$. Indeed, to prove existence of a minimizer for $S(s,\O)$, we consider a minimizing sequence $u_n$,  given by Ekeland variational principle which is bound in $H^1(\O)$ and normalized so that 
$$
\int_{\de\O}d^{-s}(\s) |u_n|^{q(s)}d\s=1.
$$
We suppose that $u_n$ converges weakly to $0$ (that is blow up occurs). Then considering the L\'{e}vy concentration function $r\mapsto \int_{\de\O\cap B_{r}(0)}d^{-s}(\s) |u_n|^{q(s)}d\s$, it easy to see, by continuity, that there exists a sequence of real number $r_n$ such that 
$$
\int_{\de\O\cap B_{r_n}(0)}d^{-s}(\s) |u_n|^{q(s)}d\s=\frac{1}{2}.
$$
Now because  $0<s$, we have $q(s)< 2^{\sharp}$ so  that by compactness $u_n\to 0$ in $L^{q(s)}(\de\O)$. Using this we show that up to a subsequence $r_n\to 0$ as $n\to \infty$. Now scaling $u_n$ with these parameters $r_n$ and making change of coordinates, we find out new sequence of functions $w_n $ for which  their mass  concentrate at a half ball centred at the origin. Namely
$$
\int_{B^N_{{r_0} } }|\ti{z}|^{-s} {w_n}^{q(s)} d\tz=\frac{1}{2} (1+O(r_n)) .
$$
 Cutting-off $w_n$ by a function $\eta_n$ near the origin and using further analysis, we see that $\eta_n w_n$
converges in $\calD^{1,2}(\ov{\R^{N+1}_+})$ to a function $w\neq0$ satisfying
\begin{align*}
 \begin{cases}
\D w=0 &\quad\textrm{ in }  \R^{N+1}_+, \vspace{3mm}\\
-\frac{\de w}{\de z^1}=S(s,\O) |\ti{z}|^{-s}|w|^{q(s)-2} w &\quad\textrm{ on }  \de \R^{N+1}_+,\vspace{3mm}
\\
\displaystyle\int_{ \de \R^{N+1}_+}|\ti{z}|^{-s}|w|^{q(s)}\, d\ti{z}\leq 1,&  \vspace{3mm},  \\
w\neq 0.
  \end{cases}
\end{align*}
This then  implies that $S(s,\O)\geq S(s)$.

\bigskip

\noindent
   \textbf{ Acknowledgments} \\
 This work is supported by the Alexander von Humboldt foundation and the    German Academic
Exchange Service (DAAD). Part of this work was done while the authors was visiting the International Center for Theoretical Physics (ICTP)
in the  Simons associateship program. The authors thank the anonymous  referee for carefully reading the first version of the manuscript and for his/her useful comments.

\section{Tool box} 
\subsection{Existence of ground states in $\R^{N+1}_+$}\label{ss:ex-gs}
We start with a proof of existence of minimizers for $S(s)$, with $s\in(0,1)$, which might be of interest in the study of Hardy-Sobolev inequalities and different from the one of \cite{Lieb}.
\begin{Theorem}\label{th:exist-s01}
 Let $s\in(0,1)$.   Then $S(s)$ has a positive minimizer $w\in \calD$ which satisfies
\begin{align*}
\begin{cases}
\displaystyle \D w =0&\qquad \textrm{ in } \R^{N+1}_+, \vspace{3mm}\\
\displaystyle -\frac{\de w}{\de z^1}=S(s) w^{q(s)-1} &\qquad \textrm{ on } \R^N, \vspace{3mm}\\
\displaystyle\int_{\R^{N}} w^{q(s)}\,dx=1.&
\end{cases}
\end{align*}

 \end{Theorem}
\proof
Recall that  $\calD:=\calD^{1,2}(\ov{\R^{N+1}_+})$.  
Define the functionals	 $\Phi,\Psi: \calD\to \R$ by
$$
\Phi(w):=\frac{1}{2} \int_{\R^{N+1}_+}|\n w|^2dx 
$$
and
$$
\Psi(w)=\frac{1}{q(s)}\int_{\de \R^{N+1}_+}|z|^{-s} |w|^{q(s)}dz.
$$
By Ekeland variational principle there exits a minimizing sequence  $w_n$  for the quotient $ S:=S(s)$ such that
\be\label{eq:wenorm}
\int_{\de \R^{N+1}_+}|z|^{-s} |w_n|^{q(s)}d z=1,
\ee
$$
\Phi(w_n)\to  \frac{1}{2}S
$$
and
\be\label{eq:weps-stf}
\Phi'(w_n)-S\Psi'(w_n)\to 0\quad\textrm{ in } \calD',
\ee
where $\calD' $ denotes the dual of  $\calD$.
We have that
\be\label{eq:webndH1}
  \int_{ \R^{N+1}_+}|\n w_n|^2dz  \leq C.
\ee
We define the Levi-type concentration function: for $r>0$
$$
Q(r):=\int_{B^N_r}|\tz|^{-s} |w_n|^{q(s)}d \tz.
$$
By continuity and \eqref{eq:uenorm} there exists $r_n>0$ such that
$$
Q(r_n):=\int_{ B^N_{r_n}}|\tz|^{-s} |w_n|^{q(s)}d \tz=\frac{1}{2}.
$$
Let $v_n(z)=r_n^{\frac{N-1}{2}}w_n(r_n z)$. It is easy to check  that for every $s\in [0,1]$
$$
\int_{ \R^{N+1}_+}|\n w_n|^2dz=\int_{ \R^{N+1}_+}|\n v_n|^2dz,\quad \int_{ \R^N}|\tz|^{-s} |w_n|^{q(s)}d \tz=\int_{ \R^N}|\tz|^{-s} |v_n|^{q(s)}d \tz
$$
and
\be\label{eq:Levi}
\int_{ B^N_1}|\tz|^{-s} |v_n|^{q(s)}d \tz=\frac{1}{2}.
\ee
Hence $v_n$ is a minimizing sequence. In particular  ${v_n} \rightharpoonup v$ for some $v$ in $ \calD$. 
We wish to show that $v\neq0$. If not then $v_n\to0$ in $L^2_{loc}(\R^{N+1}_+)$ and in $L^2_{loc}(\R^{N})$.
Let $\vp\in C^\infty_c({B_1})$ and $\vp\equiv 1$ on $B_{\frac{1}{2}}$.
 Using $\vp^2 v_n$ as test in  \eqref{eq:weps-stf} and using standard integration by parts
\begin{eqnarray*}
\displaystyle\int_{ \R^{N+1}_+}|\n (\vp v_n)|^2dz
 &=&
\displaystyle S(s)\int_{ \R^N}|z|^{-s} |v_n|^{q(s)-2}|\vp v_n|^2d \tz+o(1)\\
&\leq&\displaystyle\frac{S(s)}{2^{\frac{q(s)-2}{q(s)}}}
\displaystyle\left(   \int_{ \R^N}|\tz|^{-s} |\vp v_n|^{q(s)}d \tz  \right)^{ \frac{2}{q(s)}}+o(1),
\end{eqnarray*}
where we used \eqref{eq:Levi}. By \eqref{eq:CKNtrace} we deduce that
$$
S(s) \left(   \int_{ \R^N}|\tz|^{-s} |\vp v_n|^{q(s)}d \tz  \right)^{ \frac{2}{q(s)}} \leq
\frac{S(s)}{2^{\frac{q(s)-2}{q(s)}}}
\left(   \int_{ \R^N}|\tz|^{-s} |\vp v_n|^{q(s)}d \tz  \right)^{ \frac{2}{q(s)}}+o(1).
$$
Since  $s\in(0,1)$, we have $ S(s)> \frac{S}{2^{\frac{q(s)-2}{q(s)}}}$ so that
$$
o(1)=\int_{ \R^N}|\tz|^{-s} |\vp v_n|^{q(s)}d \tz=\int_{ B^N_1}|\tz|^{-s} |v_n|^{q(s)}d \tz +o(1)
$$
because $q(s)<2^{\sharp}$. We are therefore in contradiction with \eqref{eq:Levi}. Therefore $v\neq0$ is a minimizer. Standard arguments show that
 $v^+=\max(v,0)$ is also a minimizer and the proof is complete by the maximum principle.
\QED
 %
\subsection{Symmetry and decay estimates of ground states}
\begin{Theorem}\label{th:sym-dec}
 Let $N\geq2$  and let $w\in \calD$  such that $w>0$ and 
\begin{align}\label{eq:wgrst}
\begin{cases}
\D w=0& \quad\textrm{ in }  \R^{N+1}_+,\vspace{3mm}\\
-\de_{z^1}w:=-\frac{\de w}{\de z^1}=S(s) |\ti{z}|^{-s}w^{q(s)-1}& \quad\textrm{ on } \de \R^{N+1}_+.
  \end{cases}
\end{align}
Then we have:\\
(i) \,$w=w(z)$ only depends on $z^1$ and $|\ti{z}|$, and $w$ is
strictly decreasing in $|\ti z|$.\\
\noindent
(ii) $w(z)\leq \frac{C}{1+|z|^{N-1}}$ for all $z\in \ov{\R^{N+1}_+},$ for some positive constant  $C$.
\end{Theorem}
%
\proof
(i) \, For simplicity, we write $S$, $q$,   instead of $S(s)$,
$q(s)$. We first show that
\begin{equation}
  \label{eq:2}
\text{$w$ is strictly decreasing in $z^{N+1}$ in $\R^{N+1}_+ \setminus
  \{z^{N+1}=0\}$.}
\end{equation}
This will be shown with a variant of the moving plane method, see
 \cite{A,S,GNN1,GNN2,BN}. 
For $\lambda >0$, we consider the reflection $\R^{N+1} \to \R^{N+1}$, $z \mapsto
z_\lambda$ at the
hyperplane $\{z \in \R^{N+1}\::\: z^{N+1}=\lambda\}$. Moreover, we let
$H^\lambda:= \{z \in {\R^{N+1}_+}\::\: z^{N+1}> \lambda\}$, and we
define
$$
u^\lambda: \overline {\R^{N+1}_+ \cap H^{\lambda}} \to \R, \qquad
u^\lambda(z)=w(z_\lambda)-w(z).
$$
Then $u^\lambda$ is harmonic in
$\R^{N+1}_+ \cap H^{\lambda}$, and it satisfies
$$
u^\lambda(z)=0 \quad \text{on $\R^{N+1}_+ \cap \partial H^\lambda$}
$$
as well as
$$
-\partial_{z^1} u^\lambda(z) =  + S\bigl(\frac{w^{q-1}(z_\lambda)}{|z_\lambda|^s}
- \frac{w^{q-1}(z)}{|z|^s} \bigr) \quad \text{on $\R^N \cap H^\lambda$.}
$$
Let $u^\lambda_-= \min \{u_\lambda,0\}$.  Then
\begin{align*}
\int_{H^\lambda} &|\nabla u^\lambda_-|^2\,dz = \int_{H^\lambda} \nabla
u^\lambda \nabla u^\lambda_-\,dz= - \int_{\R^N \cap
  H^\lambda} \partial_{z^1} u^\lambda u^\lambda_-\,d\sigma(z)\\
&=  S\int_{\R^N  \cap
  H^\lambda} u^\lambda_  + \bigl(\frac{w^{q-1}}{|z_\lambda|^s}
- \frac{w^{q-1}(z)}{|z|^s}\bigr) \,d\sigma(z) \nonumber \\
&\le S  \int_{\R^N  \cap
  H^\lambda} u^\lambda_ -
|z|^{-s}[w^{q-1} (z)-w^{q-1} (z_\lambda)] \,d\sigma(z)\\
&\le  (q-1)S \int_{\R^N  \cap  H^\lambda} |u^\lambda_-(z)|^2  
|z|^{-s} w^{q-2}(z) \,d\sigma(z).
\end{align*}
In the last step we used that if $w(z_\lambda) \le w(z)$ then 
$$
w^{q-1} (z) - w^{q-1} (z_\lambda) \le (q-1) w^{q-2}(z)
[w(z)-w(z_\lambda)] = - (q-1) u^\lambda_-(z) w^{q-2}(z) 
$$
by the convexity of the function $t \mapsto
t^{q-1}$ on $(0,\infty)$. Using H{\"o}lder's inequality, we conclude that
\begin{equation}
  \label{eq:1}
\int_{H^\lambda}  |\nabla u^\lambda_-|^2  
 \le  c(\lambda)
\Bigl(\int_{\R^N  \cap H^\lambda}|z|^{-s}|u^\lambda_-|^q\,d \sigma(z) \Bigr)^{\frac{2}{q}}
\end{equation}
with
$$
c(\lambda) =  (q-1)S\Bigl(
  \int_{M_\lambda}|z|^{-s}w^{q}(z)\,d\sigma(z) \Bigr)^{\frac{q-2}{q}}\quad \text{and}\quad M_\lambda:= \{z \in \R^N  \cap
  H^\lambda\::\: u(z) > u(z_\lambda)\} 
$$
for $\lambda >0$. Since $c(\lambda) \to 0$ as $\lambda \to \infty$, we have $c(\lambda) <
S$ and therefore $u^\lambda_- \equiv 0$ in $H_\lambda \cap \R^{N+1}_+$
for $\lambda>0$ sufficiently large. As a consequence,
$$
\lambda^*:= \inf \{\lambda >0 \::\: \text{$w(z) \le
  w(z_{\lambda'})$ for all $z \in H^{\lambda'} \cap \R^{N+1}_+$ and
  all $\lambda' \ge \lambda$}\} < \infty.
$$
We claim that $\lambda^*= 0$. Indeed, if, by contradiction,
$\lambda^*>0$, then $u^{\lambda^*}$ is a nonnegative harmonic function
in $\R^{N+1}_+ \cap H^{\lambda^*}$ satisfying $u^{\lambda^*} =0$ on
$\R^{N+1}_+ \cap \partial H^\lambda$ and
$$
-\partial_{z^1} u^{\lambda^*}(z)=   \frac{w^{q-1}(z_{\lambda^*})}{|z_{\lambda^*}|^s}
- \frac{w^{q-1}(z)}{|z|^s}  \quad \text{on $\R^N \cap H^{\lambda^*}$,}
$$
where the last quantity is strictly positive whenever
$w(z_{\lambda*})>0$. Consequently, unless $w \equiv 0$,
$u^{\lambda^*}$ must be strictly positive in $\overline{\R^{N+1}_+}
\cap H^{\lambda^*}$ by the strong maximum principle.
We then choose a sufficiently large set $D$ compactly contained in
$\R^{N} \cap H^{\lambda^*}$ such that
$$
(q-1)S\Bigl(
  \int_{\R^{N} \cap H^{\lambda^*} \setminus
      D}|z|^{-s}w^{q}(z)\,d\sigma(z) \Bigr)^{\frac{q-2}{q}} < S.
$$
Then, for $\lambda<\lambda^*$ close to $\lambda^*$,  we have
$D  \subset \R^{N} \cap H^{\lambda}$,
$$
(q-1)S\Bigl(
  \int_{\R^{N} \cap H^{\lambda} \setminus
      D}|z|^{-s}w^{q}(z)\,d\sigma(z) \Bigr)^{\frac{q-2}{q}} < S.
$$
and $u^\lambda >0$ in $D$. As a consequence, $c(\lambda)< S$
for $\lambda<\lambda^*$ close to $\lambda^*$ because $M_\l\subset\R^{N} \cap H^{\lambda} \setminus
      D $. By (\ref{eq:1}) we have 
$u^\lambda \ge 0$ in $H^{\lambda} \cap \R^{N+1}_+$ for
$\lambda<\lambda^*$ close to $\lambda^*$,
 contrary to the definition of $\lambda^*$. We therefore conclude that
 $\lambda^*=0$, and this shows (\ref{eq:2}).\\
Repeating the same argument for the functions $z \mapsto w(z^1,A
\tilde {z})$, where
$A \in O(N)$ is an $N$-dimensional rotation, we conclude that $w$ only
depends on $z^1$ and $|\ti{z}|$, and $w$ is
strictly decreasing in $|\ti z|$. This ends the proof of (i).\\
%
To prove (ii), we write the \eqref{eq:wgrst} as
\begin{align*}
\begin{cases}
\D w =0& \quad \textrm{ in } \R^{N+1}\\
-\de_{z^1}w= a(z) w & \quad \textrm{ on } \R^N .
\end{cases}
\end{align*}
Where $a=S |z|^{-s} w^{q-2}\in L^p_{loc}(\R^N)$ for some $p>{N}$. Therefore by \cite{JLX}, we have that $w \in L^\infty_{loc} (\R^{N+1}_+)$. Now since \eqref{eq:wgrst} is invariant under Kelvin transform, we get immediately the result.
\QED
%
%
\subsection{Geometric preliminaries}
We let   $E_i$, ${i=2,\dots,N+1}$ be an  orthonormal basis of $ T_0\de\O$, the tangent plane of $\de \O$ at 0. 
 We will consider the Riemaninan manifold
 $(\de\O,\ti{g})$
where $\ti{g}$ is the Riemannian metric  induced by $\R^{N+1} $ on $\de \O$.
We  first introduce  geodesic normal coordinates in a neighborhood (in $\de\O$) of
$0$ with coordinates $ y'=(y^2,\dots,y^{N+1})\in\R^N.$ We set
$$
f(y'):=\textrm{Exp}_{0}^{\de\O}\left(\sum_{i=2}^{N+1}y^iE_i\right).
$$
 It is clear that the geodesic distance $d$ satisfies
\be\label{eq:dftiy}
d({f(\ti{y})})=|\ti{y}|.
\ee
In addition the above  choice of coordinates induces coordinate vector-fields on $\de\O$:
 $$
Y_i(y')=f_{*}(\partial_{y^i}),\quad \textrm{ for } i=2,\dots,N+1.
$$
Let $\ti{g}_{ij}=\la Y_i,Y_j\ra $, for $i,j=2,\dots,N+1$, be the component of the metric $\ti{g}$.
We have near the origin 
$$
\ti{g}_{ij}=\delta_{ij}+O(|y|^2).
$$
We denote by $ N_{\de\O}$ the  unit normal vector field along $\de\O$ interior to $\O$. Up to rotations,
 we will assume that $N_{\de\O}(0)=E_1 $.
 For any vector field $Y$ on $T\de\O$, we define $H(Y)=d N_{\de\O}[Y]$.
 The mean curvature  of $\de\O$ at 0 is given by 
$$
H_{_{\de\O}}(0)=\sum_{i=2}^{N+1}\la H(E_i), E_i\ra.
$$
 Now consider a local parametrization of a neighbourhood of $0$ in $ \R^{N+1}$ defined as
$$
F(y):={f(\ti{y})}+y^{1}N_{\de\O}(f(\ti{y})),\qquad y=(y^1,\ti{y})\in B_{r_0},
$$
where $B_{r_0} $ is a small ball centred at 0.
This yields the coordinate vector-fields in $\R^{N+1} $,
\begin{eqnarray*}
Y_i(y)&:=&F_{*}(\partial_{y^i}) \qquad i=1,\dots,N+1.
 \end{eqnarray*}
Let  $g_{ij}=\la Y_i,Y_j\ra $, for $i,j=1,\dots,N+1$, be the component of the flat metric $g$. It follows that 
$$
g_{ij}=\ti{g}_{ij}+ 2\langle H(Y_i),Y_j\rangle y^{1}+O(|y|^2).
$$
We have the following expansion of the metric. See for instance   \cite{MMF-PJM} for the proof.
\begin{Lemma}\label{lem:met-exp}
In a small ball  $B_{r_0}$ centered at $0$, 
\begin{eqnarray*}
g_{ij}&=& \delta_{ij}+2\langle H(E_i),E_j\rangle y^{1}+ O(|y|^2);\\[3mm]
g_{i1}&=& 0;\\[3mm]
g_{11}&=&1.
\end{eqnarray*}
\end{Lemma}
%
%
We now prove that Hardy ($s=1$) and Hardy-Sobolev trace inequality with singularity at the boundary and involving the geodesic boundary distance function  hold.
\begin{Lemma}\label{lem:HSt}
Suppose that $s\in[0,1]$. Then there exists a positive constant $C(s,\O)$ such that 
We have
$$
C(s,\O){\left(\int_{\de\O}d^{-s}(\s) |u|^{q(s)}d\s\right)^{\frac{2}{q(s)}}} \leq \int_{\O}|\n u|^2dx
+\int_{\O}| u|^2dx \qquad \forall u\in H^1(\O).
$$
\end{Lemma}
\proof
Let $u\in H^ 1(\O)$ and  pick $\eta\in C^\infty_c(F(B_{r_0}))$  such that $\eta\equiv 1$ on $F(B_{\frac{r_0}{2}})$
 and $\eta\equiv 0 $  on $F(B_{{r_0}})$. Then, by using the Sobolev trace inequality, \eqref{eq:CKNtrace} and Young's inequality, we get 
\begin{align*}
{\left(\int_{\de\O}d^{-s}(\s) |u|^{q(s)}d\s\right)^{\frac{2}{q(s)}}}&= {\left(\int_{\de\O}d^{-s}(\s) |\eta u+(1-\eta)u|^{q(s)}d\s\right)^{\frac{2}{q(s)}}}\\
&\leq 2 {\left(\int_{\de\O}d^{-s}(\s) |\eta u|^{q(s)}d\s\right)^{\frac{2}{q(s)}}}+C {\left(\int_{\de\O}  | u|^{q(s)}d\s\right)^{\frac{2}{q(s)}}}\\
&\leq C\left( \int_{\R^{N}}|y|^{-s}|(\eta u)(F(y))|^{q(s)}\,dy\right)^{\frac{2}{q(s)}} + C \|u\|^2_{H^1(\O)}\\
&\leq  C \frac{1}{S(s)} \int_{\R^{N+1}_+}|\n_y(\eta u)(F(y))|^2\,dy  + C \|u\|^2_{H^1(\O)}.
\end{align*}
Now by Lemma \ref{lem:met-exp} and some integration by parts, we deduce that 
$$
 \int_{\R^{N+1}_+}|\n_y(\eta u)(F(y))|^2\,dy\leq  C \|u\|^2_{H^1(\O)}
$$
and the proof is complete.
\QED
\subsection{Comparing  $ S(s,\O)$ and $ S(s)$}
\begin{Lemma}\label{lem:expS}
Let $\O\subset\R^{N+1}$ be a Lipschitz domain which is smooth at $0\in\de\O$.
We have the  following expansion
 \begin{eqnarray*}
S(s,\O)
&\leq & S(s)+\e C_{N,\e} H_{\de\O}(0)+ O(\rho(\e)),
 \end{eqnarray*}
where
$$
C_{N,\e}=\frac{N-2}{N}\displaystyle\int_{ B_{\frac{r_0}{\e}}^+  }z^1|{\n}_{\ti{z}}w|^2dz+
\displaystyle\int_{B_{\frac{r_0}{\e}}^+ }z^1|\de_{z^1}w|^2dz
$$
and
\begin{eqnarray*}
 \rho(\e)&=& \displaystyle\e^2 \int_{ B^+_{\frac{r_0}{\e}}}|z|^2|\n w|^2dz+\e^2 \int_{ \frac{r_0}{2\e}<|z|<\frac{r_0}{\e}} w^2dz
+\e \int_{ \frac{r_0}{2\e}<|\tz|<\frac{r_0}{\e}} w^2d\tz    \\
&+&        \e^2\displaystyle\int_{\de' B_{\frac{r_0}{\e}}^+ }|\ti{z}|^{2-s}w^{q(s)}d\ti{z}+
\int_{\de \R^{N+1}_+\setminus B_{\frac{r_0}{\e}}}|\ti{z}|^{-s}w^{q(s)}d\ti{z}   
+\e^2 \int_{ \R^{N+1}_+\cap B_{\frac{r_0}{\e}} } w^2 d{z},
\end{eqnarray*}
where $\de'B^+_r=\de B^+_r\cap \de \R^{N+1}_+$.
\end{Lemma}
\proof
 Let   $w\in\calD$, $w>0$  be the minimizer for $S(s)$ normalized so that $$\int_{ \de \R^{N+1}_+}|\ti{z}|^{-s}w^{q(s)}=1. $$
%
%
Define
$$
v_\e(F(y))=\e^{\frac{1-N}{2}}w\left(\frac{y}{\e}\right)\quad y\in \ov{B_{r_0}^+}.
$$
Let $\eta\in C^\infty_c(F(B_{r_0}))$  such that $\eta\equiv 1$ on $F(B_{\frac{r_0}{2}})$
 and $\eta\leq 1 $ on $\R^{N+1}$.
 We let
$$
{u_\e}(F(y))=\eta(F(y)) v_\e(F(y)).
$$

We have
\begin{eqnarray*}
 \displaystyle\int_\O|\n {u_\e}|^2d x&=& \displaystyle\int_\O\eta^2|\n v_\e|^2d x-\int_\O(\D \eta)\eta v_\e^2d x
+\displaystyle\int_{\de\O}\eta\frac{\de \eta}{\de \nu} v_\e^2d \s\\
&\leq&\displaystyle\int_{\O \cap F(B_{{r_0}})}|\n v_\e|^2d x+C\int_{\O\cap F(B_{{r_0})}\setminus F( B_{\frac{r_0}{2}}) } v_\e^2d x+
C\int_{\de\O\cap F(B_{{r_0}})\setminus F( B_{\frac{r_0}{2}}) } v_\e^2d \s \\
&=&\displaystyle \int_{  B_{\frac{r_0}{\e}}^+ }g^{ij}(\e z)w_iw_j\sqrt{|g|}(\e z)dz+O(\rho_1(\e)),
\end{eqnarray*}
where
$$
\rho_1(\e)=\e^2 \int_{ \frac{r_0}{2\e}<|z|<\frac{r_0}{\e}} w^2dz+ \e^2 \int_{ \frac{r_0}{2\e}<|\tz|<\frac{r_0}{\e}} w^2d\tz.
$$
Notice that
$$
g^{ij}(\e z)w_iw_j=|\n w|^2-2\e z^1\la H({\n}_{\ti{z}} w), {\n}_{\ti{z}} w\ra +O(\e^2|z|^2|\n w|^2)
$$
and
$$
\sqrt{|g|}(\e z)=1+\e z^1 H_{\de\O}(0) +O(\e^2|z|^2).
$$
Using this with the fact that $w(z)=w(z^1,|\ti{z}|)$, we get
\begin{eqnarray*}
\displaystyle\int_{ B_{\frac{r_0}{\e}}^+}g^{ij}(\e z)w_iw_j\sqrt{|g|}(\e z)dz
&\leq &\int_{  B_{\frac{r_0}{\e}}^+}|\n w|^2dz+\e C_{N,\e} H_{\de\O}(0)+O(\rho_2(\e)),
\end{eqnarray*}
where
$$
\rho_2(\e)=\e^2 \int_{ B_{\frac{r_0}{\e}}^+}|z|^2|\n w|^2dz.
$$
Hence we obtain
$$
 \displaystyle\int_\O|\n {u_\e}|^2d x\leq \int_{\R^{N+1}_+}|\n w|^2\,dz+\e C_{N,\e} H_{\de\O}(0)+O(\rho_2(\e))+O(\rho_1(\e)).
$$
On the other hand by  \eqref{eq:dftiy}, we have 
\begin{eqnarray*}
 \displaystyle\int_{\de\O}d^{-s}(\s){u_\e}^{q(s)}&=&
\int_{\de\R^{N+1}_+}d\left(\frac{F(0,\e\ti{z})}{\e}\right)^{-s}\eta^{q(s)}(\e \ti{z})w^{q(s)}\sqrt{|g|}(0,\e\ti{z})d\ti{z}\\
&=&\int_{\de\R^{N+1}_+}\left|\ti{z}\right|^{-s}(\e \ti{z})w^{q(s)}\sqrt{|g|}(0,\e\ti{z})d\ti{z}\\
&\,\,\,&-\int_{\de\R^{N+1}_+}\left|\ti{z}\right|^{-s}(1-\eta^{q(s)})(\e \ti{z})w^{q(s)}\sqrt{|g|}(0,\e\ti{z})d\ti{z}\\
&=&\int_{ \de \R^{N+1}_+}|\ti{z}|^{-s}w^{q(s)}d\ti{z}+O(\rho_3(\e)),
\end{eqnarray*}
where
\begin{eqnarray*}
\rho_3(\e)&=&\e^2\displaystyle\int_{\de' B_{\frac{r_0}{\e}}^+}|\ti{z}|^{2-s}w^{q(s)}d\ti{z}+
\int_{\de \R^{N+1}_+\setminus B_{\frac{r_0}{\e}}}|\ti{z}|^{-s}w^{q(s)}d\ti{z}.
\end{eqnarray*}

The lemma then follows by putting $\rho(\e)=\rho_1(\e)+\rho_2(\e)+\rho_3(\e)$.
\QED
\begin{Proposition}\label{prop:stric}
Let $\O\subset\R^{N+1}$ be a Lipschitz domain which is smooth at $0\in\de\O$.
 Suppose that $N\geq3$ and $s\in[0,1)$. Assume that    $H_{\de\O}(0)<0$.
Then   $S(s,\O)< S(s)$.
\end{Proposition}
\proof
Consider $w\in \calD$ given by Theorem \ref{th:exist-s01} the positive  minimizer for $S(s)$. By  Theorem \ref{th:sym-dec} we have that  $w(z)=\o(z^1,|\tz|)$  and  
\begin{align}\label{eq:wgrstcc}
 \begin{cases}
\D w=0& \quad\textrm{ in }  \R^{N+1}_+\\
-\frac{\de w}{\de z^1}=S(s) |\ti{z}|^{-s}w^{q(s)-1} &\quad\textrm{ on } \de \R^{N+1}_+
\\
\displaystyle\int_{ \de \R^{N+1}_+}|\ti{z}|^{-s}w^{q(s)}=1.&
 \end{cases}
\end{align}
In addition, thanks to Theorem \ref{th:sym-dec},  we  have 
\be\label{eq:decw}
w(z)\leq \frac{C}{1+|z|^{N-1}} \qquad \textrm{  for all $z\in\R^{N+1}_+$.}
\ee 
For $s=0$, we consider  the Escobar-Beckner (see \cite{Esc}, \cite{Beck}) function
\be\label{eq:EB}
w(z):=c_n\frac{1}{\left(1+|z|^2\right)^{\frac{N-1}{2}}},
\ee
with $c_n=\frac{2}{N-1}|S^N|^{\frac{-1}{N}}$,  which uniquely minimizes $S(0)$ up to translations.
\\
Let $\vp$ be a nonnegative radially symmetric cut-off function in $\R^{N+1} $ such that $\vp\leq 1$ in $\R^{N+1}_+$, $\vp\equiv 1$ on $B_{{2r_0}} $,
$\vp\equiv 0$ on $B_{{3r_0}} $ and $|\n \vp|+|\D\vp|\leq C$. Define $\vp_\e(z)=\vp(\e z)$ for all $z\in \R^{N+1}_+$. We multiply \eqref{eq:wgrstcc} by $|z|w\vp_\e$
and integrate by part to get
$$
\int_{B_{\frac{3r_0}{\e}}^+ }\vp_\e|z||\n w|^2= \int_{\de' B_{\frac{3r_0}{\e}}^+ }\vp_\e| w|^2+ \int_{\de' B_{\frac{3r_0}{\e}}^+ }\vp_\e|\tz|^{1-s}| w|^{q(s)}
+\frac{1}{2} \int_{ B_{\frac{3r_0}{\e}}^+ }w^2 \D(|z|\vp_\e).
$$
 By \eqref{eq:decw}, provided $N\geq 3$ we have
\be\label{eq:estwmz}
\int_{\de' B_{\frac{3r_0}{\e}}^+ }\vp_\e| w|^2+ \int_{\de' B_{\frac{3r_0}{\e}}^+ }\vp_\e|\tz|^{1-s}| w|^{q(s)}
\approx C+\e^{N-2}.
\ee
We also have
$$
\int_{ B_{\frac{3r_0}{\e}}^+ }w^2 \D(|z|\vp_\e)\leq C\e^2 \int_{ {\frac{2r_0}{\e} <|z|< \frac{3r_0}{\e}}  }\vp_\e| w|^2
+C\e \int_{ {\frac{2r_0}{\e} <|z|< \frac{3r_0}{\e}}  }\vp_\e| w|^2+C \int_{ B_{\frac{3r_0}{\e}}^+ }|z|^{-1}| w|^2
$$
and thus
$$
\int_{ B_{\frac{3r_0}{\e}}^+ }w^2 \D(|z|\vp_\e)
\approx C+\e^{N-2}.
$$
Using this and  \eqref{eq:estwmz} we  deduce that 
\be
\int_{B_{\frac{r_0}{\e}}^+ }|z||\n w|^2
\approx C+\e^{N-2}.
\ee
By using similar arguments as above (multiplying \eqref{eq:wgrstcc} by $|z|^2w\vp_\e$ and integrating by parts) we have
\begin{eqnarray*}
  \int_{B_{\frac{3 r_0}{\e}}^+ }|z|^2|\n w|^2&=& \int_{\de' B_{\frac{3r_0}{\e}}^+ }\vp_\e |\tz|| w|^2+ \int_{\de' B_{\frac{3r_0}{\e}}^+ }\vp_\e|\tz|^{2-s}| w|^{q(s)}
+\frac{1}{2} \int_{ B_{\frac{3r_0}{\e}}^+ }w^2 \D(|z|^2\vp_\e)\\
&\leq&\int_{\de' B_{\frac{3r_0}{\e}}^+ } |\tz|| w|^2+ \int_{\de' B_{\frac{3r_0}{\e}}^+ }|\tz|^{2-s}| w|^{q(s)}
+C \int_{ B_{\frac{3r_0}{\e}}^+ }w^2.
\end{eqnarray*}
 By \eqref{eq:decw}, the following estimates holds  
$$
\int_{\de' B_{\frac{3r_0}{\e}}^+ } |\tz|| w|^2 d\tz  \approx \int_{ \R^{N+1}_+\cap B_{\frac{r_0}{\e}} } w^2 d{z}\approx C+
\begin{cases}
\e^{N-3}, \quad N>3\\
|\log \e |, \quad N=3.
\end{cases}
$$
Now provided $N\geq 3$, we have 
$$
\displaystyle\int_{\de' B_{\frac{r_0}{\e}}^+ }|\ti{z}|^{2-s}w^{q(s)}d\ti{z} \approx
 C+
\begin{cases}
\e^{N-2-s}, \quad s\in(0,1)\\
|\log \e |, \quad s=0.
\end{cases}   
$$
We then  deduce that
$$
\int_{B_{\frac{r_0}{\e}}^+ }|z|^2|\n w|^2 \approx C+
\begin{cases}
\e^{N-3}, \quad N>3 \textrm{ and } s\in(0,1)\\
|\log \e |, \quad N=3  \textrm{ or } s=0 .
\end{cases}
$$
In addition, we have 
$$
 \int_{ \frac{r_0}{2\e}<|\tz|<\frac{r_0}{\e}} w^2d\tz  \approx\e^{N-2},
$$
$$
\int_{ \frac{r_0}{2\e}<|z|<\frac{r_0}{\e}} w^2dz
 \approx C+
\begin{cases}
\e^{N-3}, \quad N>3\\
|\log \e |, \quad N=3
\end{cases}
$$
and 
$$
\displaystyle\int_{\de\R^{N+1}_+\setminus B_{\frac{r_0}{\e}} }|\ti{z}|^{-s}w^{q(s)}d\ti{z} \approx
\e^{N-s}, \quad s\in[0,1).
%
$$
Thanks to {Lemma} \ref{lem:expS},  and the above estimates we conclude that, provided $N\geq 4$ and $s\in[0,1)$,
$$
S(s,\O)\leq S(s,0)+C_1 \e H_{\de\O}(0)+O(\e^{2})
$$
and if $ N=3$ or $s=0$, we get
$$
S(s,\O)\leq S(s)+C_1 \e H_{\de\O}(0)+O(\e^{2} |\log \e|),
$$
with $C_1 >0$.
\QED
%
%
\section{Existence of minimizer for $S(s,\O)$}
It is  clear from Proposition \ref{prop:stric} that the proof of Theorem \ref{th:exts-Om-i} is finalized by the   two results in this section.
However, we should emphasize that the argument following below works also for the pure Hardy case: $s=1$.
\begin{Proposition}\label{prop:exist_smO}
Let $\O\subset\R^{N+1}$ be a Lipschitz domain which is smooth at $0\in\de\O$.  
Let $s\in(0,1]$ and $N\geq 2$. Assume that $S(s,\O)< S(s) $. Then there exists a minimizer for $S(s,\O)$.
\end{Proposition}
\proof
We define $\Phi,\Psi: H^1(\O)\to \R$ by
$$
\Phi(u):=\frac{1}{2}\left(  \int_{\O}|\n u|^2dx 
+\int_{\O}| u|^2dx \right)
$$
and
$$
\Psi(u)=\frac{1}{q(s)}\int_{\de\O}d^{-s}(\s) |u|^{q(s)}d\s.
$$
By Ekeland variational principle there exits a minimizing sequence  $u_n$  for the quotient $ S(s,\O)=S(s,\O)$ such that
\be\label{eq:uenorm}
\int_{\de\O}d^{-s}(\s) |u_n|^{q(s)}d\s=1,
\ee
\be
\Phi(u_n)\to \frac{1}{2} S(s,\O)
\ee
and
\be\label{eq:ueps-stf}
\Phi'(u_n)-S(s,\O)\Psi'(u_n)\to 0\quad\textrm{ in } (H^{1}(\O))',
\ee
with $(H^{1}(\O))' $ denotes the dual of  $H^{1}(\O)$.
We have that
\be\label{eq:uebndH1}
 \int_{\O}|\n u_n|^2dx  +\int_{\de\O}| {u_n}|^2d\s  \leq Const.\quad \forall n\geq 1.
\ee
In particular  ${u_n} \rightharpoonup u$ for some $u$ in $ H^1(\O)$.\\
\textbf{Claim:}  $u\neq0$.\\
\noindent
Assume by contradiction that $u=0$ (that is blow up occur). By continuity, \eqref{eq:uenorm} and the fact that $s\in(0,1]$, there exits a sequence $r_{n}>0$ such that
\be\label{eq:concentr}
\int_{\de\O\cap B_{{r_n}}}d^{-s}(\s) |{u_n}|^{q(s)}d\s=\frac{1}{2}.
\ee
We now show that, up to a subsequence, $r_n\to0$. Indeed, by  \eqref{eq:uenorm} and \eqref{eq:concentr}
$$
\int_{\de\O\setminus B_{r_n}}d^{-s}(\s) |{u_n}|^{q(s)}d\s=\frac{1}{2}.
$$
 Since $q(s)<q(0)=2^\sharp$ for $s>0$, by compactness we have
$$
  r_n^s\,C\,\leq \int_{\de\O\setminus B_{r_n}} |{u_n}|^{q(s)}d\s\leq \int_{\de\O} |{u_n}|^{q(s)}d\s\to0\quad\textrm{ as } n\to \infty,
$$
for some positive constant $C$.\\
Define $ F_n(z)=\frac{1}{r_n}F(r_n z)$ for every $z\in B^+_{\frac{r_0}{r_n}} $ and put  $(g_n)_{i,j}=\la \de_i F_n, \de_j F_n \ra$. Clearly 
\be\label{eq:gntogEuc}
g_n\to g_{Euc}\quad C^1(K)\quad \textrm{ for every compact set } K\subset\R^{N+1},
\ee
where $g_{Euc}$ denotes the Euclidean metric. 
Let
$$
{w_n}(z)=r_n^{\frac{N-1}{2}}{u_n}(F(r_n z))\quad \forall z\in B^+_{\frac{r_0}{r_n}}.
$$
Then we get
$$
\int_{B^N_{{r_0} } }|\ti{z}|^{-s} {w_n}^{q(s)} d\tz=(1+o(1)) \int_{B^N_{{r_0}} }|\ti{z}|^{-s} {w_n}^{q(s)} \sqrt{|g_n|} d\tz
.
$$
Hence by \eqref{eq:concentr} we have
\be\label{eq:wenolmrqe}
\int_{B^N_{{r_0} } }|\ti{z}|^{-s} {w_n}^{q(s)} d\tz=\frac{1}{2} (1+c r_n) .
\ee
Let $\eta\in C^\infty_c(F(B_{r_0}))$, $\eta\equiv 1$ on $ F(B_{\frac{r_0}{2}})$ and  $\eta\equiv 0$ on $\R^{N+1} \setminus F(B_{{r_0}})$.
We define
$$
\eta_n(z)= \eta(F(r_nz))\quad\forall z\in \R^{N+1}.
$$
We have that
\be\label{eq:nwHbd}
\|\eta_n w_n \|_{\calD}\leq C \quad\forall n\in\N,
\ee
where as usual $\calD=\calD^{1,2}(\ov{\R^{N+1}}) $. Therefore
$$
\eta_n w_n\rightharpoonup w \quad\textrm{ in } \calD.
$$
We first show that $w\neq 0$. Assume by contradiction that $w\equiv 0$. Thus  
$w_n\to 0 $ in $L^p_{loc}( \R^{N+1}_+ )$ and in $L^p_{loc}( \de \R^{N+1}_+ )$  for every $1\leq p<2^\sharp$.
Let $\vp\in C^\infty_c(B_{\frac{r_0}{2}})$ be a cut-off function such that $\vp\equiv1$ on $B_{\frac{r_0}{4}} $ and $\vp\leq 1$ in $\R^{N+1}$.
Define $$\vp_n(F(y))=\vp(r_n^{-1}y).$$
We multiply \eqref{eq:ueps-stf} by $\vp_n^2 {u_n}$ (which is bounded in $H^1(\O)$)
 and integrate by parts to get
$$
\begin{array}{c}
\displaystyle\int_{ \O}\n u_n\n (\vp_n^2 {u_n}) dx  =
\displaystyle  S(s,\O) \int_{\de \O}d^{-s}(\s) | \vp_n {u_n}|^{q(s)-2}(\vp_n u_n)^2d\s+ o(1)\\
\hspace{5cm}\displaystyle \leq S(s,\O) \left(\int_{\de \O}d^{-s}(\s) | \vp_n {u_n}|^{q(s)}d\s\right)^{\frac{2}{q(s)}}+ o(1),
\end{array}
$$
where we have used \eqref{eq:uenorm}.
In the  coordinate system and after integration by parts, the above becomes
$$
\begin{array}{c}
\displaystyle\int_{ \R^{N+1}_+}|\n(\vp {w_n})|^2_{g_n}\sqrt{|g_n|} dz  =\displaystyle S(s,\O)\left( \int_{ \de\R^{N+1}_+}|\ti{z}|^{-s} | \vp {w_n}|^{q(s)}\sqrt{|g_n|} d\tz\right)^{\frac{2}{q(s)}}+o(1).
\end{array}
$$
Therefore, by \eqref{eq:gntogEuc}, for some constant $c>0$, we have
\be\label{eq:almsCont}
\begin{array}{c}
\displaystyle(1-c r_n)\int_{ \R^{N+1}_+}|\n(\vp {w_n})|^2 dz=\displaystyle S(s,\O)\left( \int_{\de \R^{N+1}_+}|\ti{z}|^{-s} | \vp {w_n}|^{q(s)} d\tz\right)^{\frac{2}{q(s)}}+o(1).
\end{array}
\ee
Hence by the Hardy-Sobolev trace  inequality  \eqref{eq:CKNtrace}, we get 
\be\label{eq:SalmSe}
\begin{array}{c}
\displaystyle(1-c r_n)S(s)\left(\int_{ \de\R^{N+1}_+}|\ti{z}|^{-s} | \vp {w_n}|^{q(s)}d\tz\right)^{\frac{2}{q(s)}}\hspace{4cm}\\
 \displaystyle\leq
 S(s,\O)\left(\int_{ \de\R^{N+1}_+}|\ti{z}|^{-s} | \vp {w_n}|^{q(s)}
 d\tz\right)^{\frac{2}{q(s)}}
+o(1).
\end{array}
\ee
Since  $S(s)> S(s,\O) $, we conclude that
$$
o(1)=\int_{\de \R^{N+1}_+}|\ti{z}|^{-s} | \vp {w_n}|^{q(s)}d\tz
=\int_{ B^N_{{r_0}}}|\ti{z}|^{-s} |  {w_n}|^{q(s)}d\tz+o(1)
$$
because by assumption $q(s)< 2^\sharp$. This is clearly in  contradiction with \eqref{eq:wenolmrqe} thus $w\neq0$.\\
Now pick $\phi\in C^\infty_c( \R^{N+1}\setminus \{0\} )$,
and put $ \phi_n(F(y))= \phi( r_n^{-1} y)$ for every $y\in B_{r_0}$. For $n$ sufficiently large,
 $ \phi_n\in  C^\infty_c(\ov{\O})$ and it  is bounded in $H^1(\O)$.
  We multiply \eqref{eq:ueps-stf} by $\phi_n$
 and integrate by parts to get
$$
\begin{array}{c}
\displaystyle\int_{ \O}\n u_n\n \phi_n dx =
\displaystyle  S(s,\O) \int_{\de \O}d^{-s}(\s) | {u_n}|^{q(s)-2}u_n \phi_n d\s+ o(1).
\end{array}
$$
Hence
$$
\begin{array}{c}
\displaystyle\int_{\R^{N+1}_+  }\la\n w_n,\n \phi\ra_{g_n} \sqrt{|g_n|}dz =
\displaystyle  S(s,\O) \int_{\de \R^{N+1}_+}|\tz|^{-s} | {w_n}|^{q(s)-2} w_n\phi\sqrt{|g_n|} d\tz+ o(1).
\end{array}
$$
Since $\eta_n\equiv 1$ on $B_{\frac{r_0}{2r_n}}$ and the support of $\phi$ is contained in an annulus, for $n$ sufficiently large
\begin{align*}
\displaystyle\int_{\R^{N+1}_+  }\la\n (\eta_n w_n),\n &\phi\ra_{g_n} \sqrt{|g_n|}dz \\
&=
\displaystyle  S(s,\O) \int_{\de \R^{N+1}_+}|\tz|^{-s} | {\eta_nw_n}|^{q(s)-2}\eta_n w_n \phi\sqrt{|g_n|} d\tz+ o(1).
\end{align*}
Since also $g_n$ converges smoothly to the Euclidean metric on the support of $\phi$,
 by passing to the limit, we infer that, for all $ \phi\in C^\infty_c( \R^{N+1}\setminus \{0\} )$
\be\label{eq: wstfAnn}
\displaystyle\int_{\R^{N+1}_+  }\n  w\n \phi\, dz 
=  S(s,\O) \int_{\de \R^{N+1}_+}|\tz|^{-s} | {w}|^{q(s)-2} w \phi\, d\tz.
\ee
Notice that $ C^\infty_c( \R^{N+1}\setminus \{0\}) $ is dense in $  C^\infty_c( \R^{N+1})$
with respect to the $H^{1}( \R^{N+1})$ norm when $N\geq 2$, see e.g. \cite{Maz}.
Consequently since $w\in \calD$, it follows  that  \eqref{eq: wstfAnn} 
holds for all $ \phi\in C^\infty_c( \R^{N+1})$ by \eqref{eq:CKNtrace}. We conclude that
\begin{align*}
 \begin{cases}
\D w=0 &\quad\textrm{ in }  \R^{N+1}_+, \vspace{3mm}\\
-\frac{\de w}{\de z^1}=S(s,\O) |\ti{z}|^{-s}|w|^{q(s)-2} w &\quad\textrm{ on }  \de \R^{N+1}_+,\vspace{3mm}
\\
\displaystyle\int_{ \de \R^{N+1}_+}|\ti{z}|^{-s}|w|^{q(s)}\, d\ti{z}\leq 1,&  \vspace{3mm},  \\
w\neq 0.
  \end{cases}
\end{align*}
Multiplying this equation by $w$ and integrating by parts, leads to $S(s,\O)\geq S(s)$ 
by \eqref{eq:CKNtrace} which is a contradiction and thus $u=\lim u_n\neq0$
is a minimizer for $S(s,\O)$.\\
\QED

In the following we study the existence of minimizers for the Sobolev trace inequality.
\begin{Proposition}\label{prop:exist_smOs0}
Let $\O\subset\R^{N+1}$ be a Lipschitz domain which is smooth at $0\in\de\O$ and $N\geq 2$. 
Assume that $S(0,\O)< S(0) $. Then there exists a minimizer for $S(0,\O)$.
\end{Proposition}
\proof
Recall the Sobolev trace inequality, proved by Li and Zhu in  \cite{YanYanLi}: there  exists a positive constant $C=C(\O)$ such that for all 
$u\in H^1(\O)$, we have 
\be\label{eq:almstST}
\displaystyle  S(0) \left(\int_{\de\O}| u|^{2^{\sharp}}  d\s\right)^{2/2^{\sharp} }
 \leq  \int_{\O} |\n u|^2dx+ C \int_{\de \O}| u|^2d\s.
\ee

Now we let $u_n$ be a minimizing sequence for $S(0)$, normalized as $\|u_n\|_{L^ {2\sharp}(\de\O)}=1$.
We now show that $u=\lim u_n$ is not zero. Put $\th_n:= u_n-u$ so that $ \th_n  \rightharpoonup 0$ in $H^1(\O)$  and
$\th_n\to 0$ in $L^2(\O), L^2(\de \O)$.
Moreover by Brezis-Lieb Lemma \cite{BL} and recalling \eqref{eq:uenorm},   it holds that
\be\label{eq:BrLi}
1- \lim_{n\to \infty}\int_{\de\O}| \th_n|^{2^{\sharp}} d\s= \int_{\de\O}| u|^{2^{\sharp}} d\s.
\ee
By using \eqref{eq:almstST}, we have
\begin{align*}
 S(0,\O)& \left(\int_{\de\O}| u|^{2^{\sharp}} d\s\right)^{2/2^{\sharp} }
 \leq  \int_{\O} |\n  u|^2dx  + \int_{\de\O}| {u}|^2d\s\\
 &\leq  \int_{\O} |\n  u_n|^2dx+  \int_{\de\O}| {u_n}|^2d\s
 - \int_{\O} |\n  \th_n|^2dx+o(1) \leq  \int_{\O} |\n  u_n|^2dx+\int_{\de\O}| {u_n}|^2d\s \\
& -S(0 )  \left(\int_{\de\O}| \th_n|^{2^{\sharp}} d\s\right)^{2/2^{\sharp} } +o(1)\\
&\leq  S(0,\O) -  S(0)  \left(\int_{\de\O}| \th_n|^{2^{\sharp}} d\s\right)^{2/2^{\sharp} }+o(1).
\end{align*}
We take the limit as $n\to \infty$ and use \eqref{eq:BrLi} to get
\begin{eqnarray*}
 S(0,\O) \left(\int_{\de\O}| u|^{2^{\sharp}} d\s\right)^{2/2^{\sharp} } \leq S(0,\O) -
 S(0)  \left(1- \int_{\de\O}| u|^{2^{\sharp}} d\s\right)^{2/2^{\sharp} }.
\end{eqnarray*}
Thanks to the concavity of the function $t\mapsto t^{2/2^{\sharp}  }$, the above implies that $ \int_{\de\O}| u|^{2^{\sharp}} d\s\geq 1$ whenever
$S(0,\O)< S(0)$.  This completes  the proof.
\QED

    \label{References}

\end{document}